\documentclass[12pt,draftcls,onecolumn]{IEEEtran}

\usepackage{amsmath,amsfonts,amssymb,color}
\usepackage{tikz,subcaption,cite}
\usepackage{epstopdf}
\usepackage{epsfig} 
\usepackage{psfrag}

\newtheorem{theorem}{Theorem}
\newtheorem{corollary}{Corollary}
\newtheorem{proposition}{Proposition}
\newtheorem{lemma}{Lemma}
\newtheorem{remark}{Remark}
\newtheorem{definition}{Definition}
\newtheorem{example}{Example}

\DeclareMathOperator{\diag}{diag}

\definecolor{cyan}{rgb}{0.0, 1.0, 1.0}

\begin{document}
\title{\LARGE \bf  A Graph Theoretic Approach to the Robustness of  $k$-Nearest Neighbor Vehicle Platoons}

\author{Mohammad Pirani, Ehsan Hashemi, John W. Simpson-Porco,  Baris Fidan and Amir Khajepour 
\thanks{M. Pirani, E. Hashemi,  B. Fidan and A. Khajepour  are with the Department of Mechanical and Mechatronics Engineering, University of Waterloo, ON, Canada. email: mpirani, ehashemi, fidan, a.khajepour@uwaterloo.ca.
J. S. Porco is with the department of Electrical and Computer Engineering, University of Waterloo, ON, Canada. email: jwsimpson@uwaterloo.ca}}

\maketitle

\thispagestyle{empty}
\pagestyle{empty}


\begin{abstract}

We consider a graph theoretic approach to the performance and robustness of a platoon of vehicles, where each  vehicle communicates with its $k$-nearest neighbors. In particular, we quantify the platoon's stability margin, robustness to disturbances (in terms of system $\mathcal{H}_{\infty}$ norm), and maximum delay tolerance via graph-theoretic notions such as nodal degrees and (grounded) Laplacian matrix eigenvalues. Our results show that there is a trade-off between robustness to time delay and robustness to disturbances. Both first-order dynamics (reference velocity tracking) and second-order dynamics (controlling inter-vehicular distance) are analyzed in this direction. Theoretical contributions are confirmed via simulation results.

\end{abstract}

\section{Introduction}

Recent achievements in redefining system theoretic notions (such as stability and robustness) based on network properties have found various applications in real-world networked systems \cite{Olfati2}. One important application is in connected and cooperative vehicles, which will have great impacts on forming future generation of urban transportation \cite{Bamieh2}. Among different applications of connected vehicles, cooperative cruise control has attracted much attention. This method is concerned with controlling vehicles' velocities to minimize fuel consumption and  maintain prescribed inter-vehicular distances. However, one of the main challenges is to make these control policies resilient to external disturbances and to time delays in inter-vehicle communications. The aim of this paper is to present conditions for  the robustness of a generalized form of vehicle platooning (called $k$-nearest neighbor platoon) to communication disturbances and time delay. 

Much effort has been made in analyzing the robustness of vehicle platoons to  communication disturbances and among them is the well-known notion of string stability \cite{seiler2004disturbance}. String stability occurs if the transfer function from disturbance in the first vehicle in a platoon to state error in the last vehicle has a bounded frequency magnitude peak independent of the platoon size \cite{middleton2010string}. The notion of robustness was revisited later in terms of network coherence in the control theory literature \cite{bamieh2012coherence}. It is shown that in  1-D network topologies it is impractical to have large coherent platoons with only local feedbacks. Alternatively, optimal controllers are designed in \cite{lin2012optimal} to improve the coherence of a vehicle formation. The effect of time delay in vehicle communication on performance and stability in vehicle platoons was studied in  \cite{liu2001effects, xiao2009scalability}.

This paper is concerned with  robustness analysis of vehicle platoons to delay and communication disturbances under two policies. First we analyze the velocity tracking scenario, which is  applied to  cases where the vehicle fuel consumption is to be minimized \cite{van2015fuel, koller2015fuel}. Second, we analyze the network formation problem, where inter-vehicular distances are regulated to avoid collisions. In contrast with other works on robustness of vehicle platoons \cite{chien1992automatic, klinge2009time, hao2013stability, herman2015nonzero}, here we present graph theoretic robustness conditions for both of the above communication policies and analyze the effect of the number and location of the reference vehicles (leaders) on robustness of the vehicle network. More specifically, the contributions of this paper are the followings:
\begin{itemize}
\item We provide graph theoretic bounds for the system $\mathcal{H}_{\infty}$ norm of the {\it{velocity tracking scenario}}. Moreover, we propose  necessary and sufficient conditions for the value of the maximum constant time delay $\tau_{\rm max}$ for which the velocity tracking remains asymptotically stable.  

\item We provide graph theoretic bounds for the  system $\mathcal{H}_{\infty}$ norm  of the {\it{network formation problem}}. In addition, we introduce an upper bound for  $\tau_{\rm max}$ such that the network formation remains asymptotically stable. 
\end{itemize}

The paper is organized as follows. We start by introducing some required notations in Section \ref{sec:not} and in Section \ref{sec:influence} the network dynamics for both {\it{velocity tracking}} and {\it{network formation}} is provided. Section \ref{sec:eigenvaluesmall} briefly presents some graph theoretic bounds on the extreme eigenvalues of the grounded Laplacian matrix which are used in Sections \ref{sec:singleeee} and \ref{sec:doubleeee} to establish  conditions for the robustness of $k$-nearest neighbor platoons for both velocity tracking and network formation scenarios. In Section \ref{sec:Results} we present some simulation results and Section \ref{sec:conclusion} concludes the paper.

\section{Notations and Definitions}
\label{sec:not}

We denote an undirected graph (network) by  $\mathcal{G}=\{\mathcal{V},\mathcal{E}\}$,  where $\mathcal{V} = \{v_1, v_2, \ldots, v_n\}$ is a set of nodes (or vertices) and $\mathcal{E} \subset \mathcal{V}\times\mathcal{V}$ is the set of edges.  Neighbors of node $v_i \in \mathcal{V}$ are given by the set $\mathcal{N}_i = \{v_j \in \mathcal{V} \mid (v_i, v_j) \in \mathcal{E}\}$. The  adjacency matrix of the graph is given by a symmetric and binary $n \times n$  matrix $A$, where element $A_{ij}=1$ if $(v_i, v_j) \in \mathcal{E}$ and zero otherwise.  The degree of node $v_i$ is denoted by  $d_i \triangleq \sum_{j=1}^nA_{ij}$. For a given set of nodes $X \subset \mathcal{V}$, the {\it edge-boundary} (or just boundary) of the set is given by $\partial{X} \triangleq \{(v_i,v_j) \in \mathcal{E} \mid v_i \in X, v_j \in \mathcal{V}\setminus{X}\}$.    The Laplacian matrix of the graph is given by $\mathcal{L} \triangleq D - A$, where $D = \diag(d_1, d_2, \ldots, d_n)$.  The eigenvalues of the Laplacian are real and nonnegative, and are denoted by $0 = \lambda_1(\mathcal{L}) \le \lambda_2(\mathcal{L}) \le \ldots \le \lambda_n(\mathcal{L})$.   For a given subset $\mathcal{S} \subset \mathcal{V}$ of nodes (which we term {\it grounded nodes}), the {\it grounded Laplacian} induced by $\mathcal{S}$ is denoted by $\mathcal{L}_g$, and is obtained by removing the rows and columns of $\mathcal{L}$ corresponding to the nodes in $\mathcal{S}$. In this paper,  {\it grounded nodes} represent {\it reference vehicles}. For the case where the underlying network is connected and there exists at least one grounded node, the grounded Laplacian matrix $\mathcal{L}_g$ is a positive definite matrix \cite{PiraniSundaramArxiv}. For a given set $\mathcal{I}$, the number of members (cardinality) of the set is denoted by $|\mathcal{I}|$.

\section{Problem Statement}
\label{sec:influence}
\begin{figure}[h!]
\centering
\includegraphics[width=0.8\linewidth]{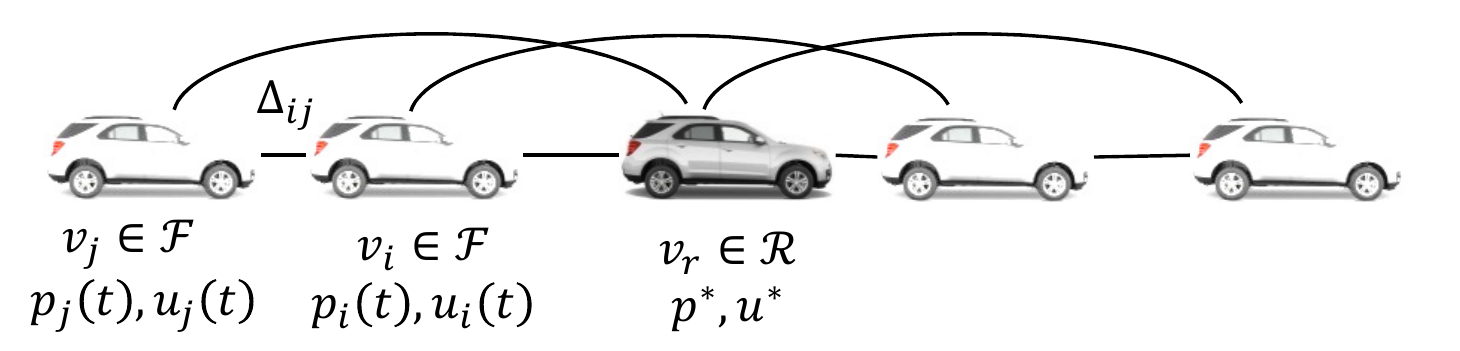}
\caption{A $2$-nearest neighbor platoon of 5 vehicles, $\mathcal{P}(5,2)$. A reference vehicle is located in the middle.}
\label{fig:faultt}
\end{figure}

Consider  a connected network of  $n$ vehicles $\mathcal{V} = \{v_1, v_2,\ldots, v_n\}$. Each vehicle $v_i \in \mathcal{V}$ is either a follower $v_i \in \mathcal{F}$ or a reference vehicle $v_i \in \mathcal{R}$. The position and longitudinal velocity of each vehicle $v_i$ is denoted by scalars $p_i$ and $u_i$, respectively, which evolve over time with  particular dynamics. In this paper $\mathcal{P}(n,k)$ denotes a platoon of $n$ vehicles  where each vehicle can communicate with its $k$ nearest neighbors from  its back and $k$ nearest neighbors from its front, for some $k \geq 1$. This is due to the limited communication range for each sensor in a vehicle and the distance between  the consecutive vehicles. An example of $\mathcal{P}(n,k)$ is shown in Fig. \ref{fig:faultt}.  

Two control objectives are addressed in this paper: (i) control the velocity of the vehicles (velocity tracking) or (ii) regulate the distance between neighboring vehicles (network formation). Each vehicle $v_i$ is governed by the second order dynamics $\ddot{p}_i(t)= q_i(t)$, or in vector notation
\begin{equation}
\ddot{\mathbf{p}}(t)=\mathbf{q}(t),
\label{eqn:law}
\end{equation}
where  $\mathbf{p}(t)=[{p}_1(t), {p}_2(t), ..., {p}_n(t)]^T$ is the vector of positions and  $\mathbf{q}(t)$ is the vector of control laws for either of the two above control objectives discussed in the following subsections.

\subsection{Velocity Tracking}
In this case, each follower vehicle  tracks a reference velocity trajectory. The desired velocity is calculated  by reference vehicles based on  minimizing fuel consumption. This yields  the following control laws for each follower and reference vehicle, \cite{Rahmani}
\begin{align}
{q}_i(t)=
  \begin{cases}
    \sum_{j\in \mathcal{N}_i}k_u(u_j(t)-u_i(t))       & \quad  \forall  v_i\in \mathcal{F},\\
   0  & \quad \forall v_i \in \mathcal{R},\\
  \end{cases}
\label{eqn:single}
\end{align}
where $k_u>0$ is the control gain. The state (velocity) of the reference vehicles (which should be tracked by the followers) is assumed to be constant and is not affected by other vehicles. 
\begin{remark}
One can define a control law ${q}_i(t)=\kappa(u^{*}-u_i(t))$ for all $v_i \in \mathcal{R}$ where $u^{*}$ is the reference velocity. For sufficiently large $\kappa$  it can be shown  (by singular perturbation analysis) that the assumption  ${q}_i(t)=0$ in \eqref{eqn:single} is valid.    
\end{remark}

 Aggregating the velocities of all followers  into a vector $\mathbf{u}_\mathcal{F}(t) \in \mathbb{R}^{|\mathcal{F}|}$, and the velocities of all  reference vehicles into a vector $\mathbf{u}_{\mathcal{R}}(t)\in \mathbb{R}^{|\mathcal{R}|}$,\footnote{Note that $\mathbf{u}_{\mathcal{R}}(t) =u^{*}\mathbf{1}_{|\mathcal{R}|\times 1}$ for all $t \ge 0$, where $u^{*}$ is a unique reference velocity. The reason of using multiple reference vehicles with the same value is to increase the network robustness as will be discussed later.} one can write the following dynamics from  \eqref{eqn:law} and \eqref{eqn:single}:
\begin{equation}
\begin{bmatrix}
      \dot{\mathbf{u}}_\mathcal{F}(t)          \\[0.3em]
       \dot{\mathbf{u}}_{\mathcal{R}}(t) 
     \end{bmatrix}=-k_u\underbrace{\begin{bmatrix}
       \mathcal{L}_g & \mathcal{L}_{12}           \\[0.3em]
       \mathbf{0} &  \mathbf{0}          
     \end{bmatrix}}_{\mathcal{L}}\begin{bmatrix}
      {\mathbf{u}}_{\mathcal{F}}(t)          \\[0.3em]
       {\mathbf{u}}_{\mathcal{R}}(t) 
     \end{bmatrix}.
\label{eqn:mat}
\end{equation}
where   $\mathcal{L}_g\in \mathbb{R}^{|\mathcal{F}|\times|\mathcal{F}|}$ is the grounded Laplacian matrix, formed by removing the rows and columns corresponding to the reference vehicles. The control law for the follower vehicles in vector form becomes
\begin{equation}
\dot{\mathbf{u}}_\mathcal{F}(t)=-k_u\mathcal{L}_g{\mathbf{u}}_{\mathcal{F}}(t)-k_u\mathcal{L}_{12}{\mathbf{u}}_{\mathcal{R}}(0).
\label{eqn:sing}
\end{equation}

\begin{remark}
 The unique steady-state solution of \eqref{eqn:sing} is ${\mathbf{u}}_{\mathcal{F}}^{ss}=-\mathcal{L}_g^{-1}\mathcal{L}_{12}{\mathbf{u}}_{\mathcal{R}}(0)$. We know that $\mathcal{L}\mathbf{1}=\mathbf{0}$ which yields $\mathcal{L}_g\mathbf{1}_{|\mathcal{F}|\times 1}+\mathcal{L}_{12}\mathbf{1}_{|\mathcal{R}|\times 1}=\mathbf{0}_{|\mathcal{F}|\times 1}$ and that results in $-\mathcal{L}_g^{-1}\mathcal{L}_{12}\mathbf{1}_{|\mathcal{R}|\times 1}=\mathbf{1}_{|\mathcal{F}|\times 1}$. Hence,   $-\mathcal{L}_g^{-1}\mathcal{L}_{12}$ is a row stochastic matrix and since there is only one reference velocity,  ${\mathbf{u}}_{\mathcal{F}}^{ss}$ attains that value. 
\end{remark}
 
 Introducing $\tilde{\mathbf{u}}_{\mathcal{F}}(t)={\mathbf{u}_{\mathcal{F}}}(t)-{\mathbf{u}}_{\mathcal{F}}^{ss}$, the error dynamics of \eqref{eqn:sing} becomes
\begin{equation}
\dot{\tilde{\mathbf{u}}}_\mathcal{F}(t)=-k_u\mathcal{L}_g{\tilde{\mathbf{u}}}_{\mathcal{F}}(t).
\label{eqn:errord}
\end{equation}

\subsection{Network Formation } 
In this case, the objective for each follower vehicle is to maintain specific distances from its neighbor vehicles.  The desired vehicle formation will be formed by a specific constant distance $\Delta_{ij}$ between vehicles $v_i$ and $v_j$, which should satisfy $\Delta_{ij}=\Delta_{ik}+\Delta_{kj}$ for every triple $\{v_i,v_j,v_k\}\in \mathcal{V}$. The control law for each follower vehicle is \cite{Hao}
\begin{equation}
q_i(t)=\sum_{j\in \mathcal{N}_i}k_p \left (p_j(t)-p_i(t)+\Delta_{ij} \right )+k_u \left (u_j(t)-u_i(t)\right),
\label{eqn:singlle}
\end{equation}
where $k_p,k_u>0$ are control gains. We define the tracking error $\tilde{p}_i(t)=p_i(t)-p^*_i(t)$, where $p^*_i(t)$ is the desired trajectory of vehicle $v_i$ which should satisfy $\Delta_{ij}=p^*_i(t)-p^*_j(t)$ for all $v_i, v_j \in \mathcal{V}$. By rewriting \eqref{eqn:singlle} we have 
\begin{equation}
\ddot{\tilde{p}}_i(t)=\sum_{j\in \mathcal{N}_i}k_p \left (\tilde{p}_j(t)-\tilde{p}_i(t) \right )+k_u \left (u_j(t)-u_i(t)\right),
\label{eqn:singaslle}
\end{equation}
which comes from the fact that the rigid formation requires $u^*_j(t)=u^*_i(t)$ which results in $\tilde{u}_j(t)-\tilde{u}_i(t)=u_j(t)-u_i(t)$. The error dynamics \eqref{eqn:singaslle} in the state space form is
\begin{equation}
\dot{\tilde{\mathbf{x}}}_{\mathcal{F}}(t)=\mathcal{B}\tilde{\mathbf{x}}_{\mathcal{F}}(t),
\label{eqn:doub}
\end{equation}
where $\tilde{\mathbf{x}}_{\mathcal{F}}=[\tilde{{p}}_1, \tilde{{p}}_2, ..., \tilde{{p}}_{|\mathcal{F}|}, \dot{\tilde{{p}}}_1, \dot{\tilde{{p}}}_2,  ..., \dot{\tilde{{p}}}_{|\mathcal{F}|}] $                                                  and   $\mathcal{B}=I_{|\mathcal{F}|\times|\mathcal{F}|}\otimes \mathcal{B}_1+\mathcal{L}_g\otimes \mathcal{B}_2$, in which $\otimes$ is Kronecker product and 
\begin{equation}
\mathcal{B}_1=\begin{bmatrix}
       0 & 1         \\[0.3em]
     0 & 0
     \end{bmatrix}, \quad \mathcal{B}_2=\begin{bmatrix}
       0 & 0         \\[0.3em]
     -k_p & -k_u
     \end{bmatrix}.
\end{equation}

\begin{remark}
Going forward we assume $k_p=k_u=1$ and we focus on the effect of the network structure (not control gains) on the robustness of vehicle platoon. The results can be easily extended for all $k_p,k_u>0$. 
\end{remark}
The following theorem, introduces the spectrum of matrix $\mathcal{B}$ in \eqref{eqn:doub} in terms of the spectrum of $\mathcal{L}_g$. 
 \begin{theorem}[\cite{Hao}]
 The spectrum of $\mathcal{B}$, $\sigma(\mathcal{B})$, is 
 \begin{align}
 \sigma(\mathcal{B})=\cup_{\substack{\\ \lambda_i\in \sigma(\mathcal{L}_g)}} \left\{\sigma \begin{bmatrix}
       0 & 1          \\[0.3em]
       -\lambda_i & -\lambda_i      
     \end{bmatrix}\right\}.
     \label{eqn:char}
 \end{align}
 Thus by forming the characteristic polynomial of \eqref{eqn:char} we have
 \begin{equation}
\lambda_i(\mathcal{B})=
\begin{cases}
  -\frac{\lambda_{i}(\mathcal{L}_g)}{2}\left(1 + \left(1-\frac{4}{\lambda_{i}(\mathcal{L}_g)}\right)^{\frac{1}{2}} \right),  & \hspace{-8 mm} 1\leq i\leq |\mathcal{F}|,\\
  -\frac{\lambda_{i-|\mathcal{F}|}(\mathcal{L}_g)}{2}\left(1 - \left(1-\frac{4}{\lambda_{i-|\mathcal{F}|}(\mathcal{L}_g)}\right)^{\frac{1}{2}} \right), & \hspace{-2 mm} i> |\mathcal{F}|,\\
  \end{cases}
\label{eqn:charact}
 \end{equation}
 where $i=1,2, ..., 2|\mathcal{F}|$. Based on the fact that $\lambda_{i-|\mathcal{F}|}(\mathcal{L}_g)$ for $|\mathcal{F}|\leq i\leq 2|\mathcal{F}|$ is the same as $\lambda_{i}(\mathcal{L}_g)$ for $1 \leq i\leq |\mathcal{F}|$,  each eigenvalue of $\mathcal{L}_g$ in \eqref{eqn:charact} forms two  eigenvalues of $\mathcal{B}$ and since $\mathcal{L}_g$ is a positive definite matrix, the real parts of all of the eigenvalues of $\mathcal{B}$ are negative.
 \label{thm:eigenthm}
 \end{theorem}

 \subsection{Robustness Notions for Vehicle Platoons}
 
\begin{figure}[h!]
\centering
\includegraphics[width=0.6\linewidth]{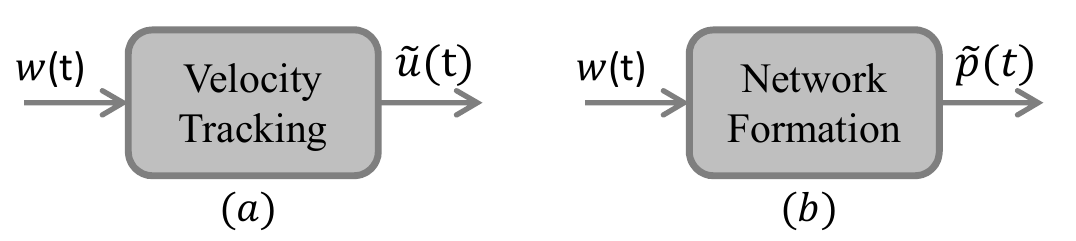}
\caption{Input and outputs of (a) error dynamics \eqref{eqn:errord} and (b) error dynamics \eqref{eqn:doub}.}
\label{fig:inputoutput}
\end{figure}

From now on we  refer to the error dynamics \eqref{eqn:errord} as {\it{velocity tracking dynamics}} and to \eqref{eqn:doub} as {\it{network formation dynamics}}. These control policies are prone to imprecisions due to  the inter-vehicle communication disturbances. Hence, both velocity tracking dynamics and network formation dynamics can be written in the following form
\begin{equation}
\dot{\mathbf{x}}(t)=\mathcal{A}\mathbf{x}(t)+\mathcal{J}\mathbf{w}(t),
\label{eqn:ladww}
\end{equation}
where $\mathbf{w}(t)$ is a vector which represents bounded disturbances. Here $\mathcal{A}=-\mathcal{L}_g, \hspace{1mm} \mathcal{J}=I_{|\mathcal{F}|\times |\mathcal{F}|}$ for the velocity tracking dynamics and $\mathcal{A}=\mathcal{B}, \hspace{1mm}\mathcal{J}=[I_{|\mathcal{F}|\times |\mathcal{F}|}  \quad \mathbf{0}_{|\mathcal{F}|\times |\mathcal{F}|}]^T$ for the network formation dynamics. As output signals of interest, we consider {\it{velocity}} for the velocity tracking dynamics and {\it{position}} for the network formation dynamics, as shown in  Fig. \ref{fig:inputoutput}. 
Based on the input-output representation of both dynamics, the robustness of vehicle platoon to disturbances is analyzed based on the $\mathcal{H}_{\infty}$ norm of the transfer function from the disturbances to the output signals. 
\begin{remark}
The notion of system $\mathcal{H}_{\infty}$ norm discussed in this paper is to address the robustness of each agent's state error (position or velocity) to external disturbances. Thus, it is different from the notion of $\mathcal{L}_2$ string stability \cite{seiler2004disturbance}, which addresses the effect of the disturbances on the first vehicle to the state error of the last vehicle in a platoon. 
\end{remark}

In addition to disturbances, the inter-vehicle communication is prone to  time delay  which may inhibit tracking or even causes instability. 
More formally, updating policies \eqref{eqn:errord} and \eqref{eqn:doub} can be in the following form
\begin{equation}
\dot{\mathbf{x}}(t)=\mathcal{A}\mathbf{x}(t-\tau),
\label{eqn:lawwdel}
\end{equation}
where $\tau \in [0,\tau_{\rm max}]$ is a constant time delay. 
\begin{remark}
For velocity tracking dynamics \eqref{eqn:errord} if each vehicle has instantaneous access to its own state, the dynamics have the form
\begin{equation}
\dot{\tilde{\mathbf{u}}}_\mathcal{F}(t)=-D_g{\tilde{\mathbf{u}}}_{\mathcal{F}}(t)+A_g{\tilde{\mathbf{u}}}_{\mathcal{F}}(t-\tau),
\label{eqn:offdiag}
\end{equation}
where $\mathcal{L}_g=D_g-A_g$. In this case  since all of the principal minors of $\mathcal{L}_g$ are nonnegative and $\mathcal{L}_g$ is non-singular, \eqref{eqn:offdiag} is asymptotically stable independent of the magnitude of the delays in the off-diagonal terms of $\mathcal{L}_g$ (Theorem 1 in \cite{Hofbauer}). Hence, for the sake of consistency, we introduce conditions for $\tau$ where \eqref{eqn:lawwdel} is stable for both dynamics \eqref{eqn:errord} and \eqref{eqn:doub}. 
\end{remark}

In the following section, a brief overview about the spectrum of the grounded Laplacian matrix is presented.



\section{Smallest and Largest Eigenvalues of $\mathcal{L}_g$}
\label{sec:eigenvaluesmall}

 Spectrum of $\mathcal{L}_g$ has a pivotal role in the performance and robustness of the velocity tracking and network formation dynamics  \eqref{eqn:errord} and \eqref{eqn:doub}.  The following theorems provide  bounds on $\lambda_1(\mathcal{L}_g)$ and $\lambda_{|\mathcal{F}|}(\mathcal{L}_g)$ based on  network properties. 

\begin{theorem}[\cite{PiraniSundaramArxiv}]
Consider a connected network $\mathcal{G}=\{\mathcal{V},\mathcal{E}\}$ with a set of reference vehicles $\mathcal{R} \subset \mathcal{V}$.  Let $\mathcal{L}_g$ be the grounded Laplacian matrix for $\mathcal{G}$.  Let $\beta_i=|\mathcal{N}_i\cap \mathcal{R}|$ be the number of reference vehicles in follower $v_i$'s neighborhood. Then
\begin{equation}
 \min_{i\in \mathcal{F}}\{\beta_i\} \leq \lambda_1(\mathcal{L}_g) \leq   \frac{|\partial \mathcal{R}|}{|\mathcal{F}|}\leq \max_{i\in \mathcal{F}}\{\beta_i\}\leq |\mathcal{R}|.
\label{eqn:maineq}
\end{equation}
\label{thm:main}
\end{theorem}
\begin{theorem}[\cite{ArxiveRobutness}]
Consider a connected graph $\mathcal{G}=\{\mathcal{V},\mathcal{E}\}$ with a set of reference vehicles $\mathcal{R} \subset \mathcal{V}$.  Let $\mathcal{L}_g$ be the grounded Laplacian matrix for $\mathcal{G}$. For the largest eigenvalue of $\mathcal{L}_g$ we have
\begin{align}
d_{\rm max}^{\mathcal{F}} \leq \lambda_{|\mathcal{F}|}(\mathcal{L}_g) \leq 2d_{\rm max}^{\mathcal{F}},
\label{eqn:up}  
\end{align}
where $d_{\rm max}^{\mathcal{F}}=\max_{v_i\in \mathcal{F}}d_i$, is the maximum degree over the follower vehicles.
\label{thm:bounddmax}
\end{theorem}

\section{Robustness of Velocity Tracking Dynamics}
\label{sec:singleeee}
In the previous section, some useful spectral properties of the grounded Laplacian matrix $\mathcal{L}_g$ were introduced. In this section, we use those results to give graph theoretic conditions for the stability margin of dynamics \eqref{eqn:errord} and its robustness to disturbances and time delay. 

\subsection{Stability Margin and Robustness to Disturbances}

The stability margin  of  \eqref{eqn:errord} is determined by $\lambda_1(\mathcal{L}_g)$. Hence, the graph theoretic bounds provided in \eqref{eqn:maineq} can be considered as bounds on the stability margin accordingly. Now suppose that there exists an additional disturbance vector $\mathbf{w}(t)$ in velocity tracking dynamics in the form of \eqref{eqn:ladww} which is 
\begin{equation}
\dot{\tilde{\mathbf{u}}}_\mathcal{F}(t)=-\mathcal{L}_g{\tilde{\mathbf{u}}}_{\mathcal{F}}(t)+\mathbf{w}(t).
\label{eqn:errodrd}
\end{equation}
 The transfer function of \eqref{eqn:errodrd} is $G(s)= (sI+\mathcal{L}_g)^{-1}$. Here, the system  $\mathcal{H}_{\infty}$ norm of  \eqref{eqn:errodrd} is considered, defined as $||G||_{\infty} \triangleq \sup_{\omega\in \mathbb{R}}{\lambda_{max}^{\frac{1}{2}}(G^*(j\omega)G(j\omega))}$ \cite{Doyle}, as introduced in the following proposition. 
\begin{proposition}[\cite{piranicdc}]
The system  $\mathcal{H}_{\infty}$ norm  of  \eqref{eqn:errodrd} is $||G||_{\infty}=\frac{1}{\lambda_1({L}_g)}$. 
\label{prop:hinferds}
\end{proposition}

Based on Proposition \ref{prop:hinferds} and Theorem \ref{thm:main}, the following bounds for  $\mathcal{H}_{\infty}$ norm of \eqref{eqn:errodrd} can be written:
\begin{equation}
 \frac{1}{\max_{i\in \mathcal{F}}\{\beta_i\}} \leq \frac{|\mathcal{F}|}{|\partial \mathcal{R}|} \leq ||G||_{\infty} \leq \frac{1}{\min_{i\in \mathcal{F}}\{\beta_i\}}.
\label{eqn:hinfff}
\end{equation}
For the case where $\min_{i\in \mathcal{F}}\{\beta_i\}=0$, the upper bound in \eqref{eqn:hinfff} is infinity. According to \eqref{eqn:hinfff}, we have the following corollary. 
\begin{corollary}
Consider a vehicle platoon $\mathcal{P}(n,k)$ with reference vehicle set $\mathcal{R}$ and follower set $\mathcal{F}$. Necessary and sufficient conditions for $\mathcal{P}(n,k)$ to have  $||G||_{\infty}< \gamma$ are to have $\max_{i\in \mathcal{F}}\{\beta_i\}> \lfloor \frac{1}{\gamma} \rfloor$ and $\min_{i\in \mathcal{F}}\{\beta_i\}> \lceil \frac{1}{\gamma} \rceil$, respectively. 
\label{cor:corolllll}
\end{corollary}

The following theorem  addresses  necessary and sufficient conditions for the number of reference vehicles in $\mathcal{P}(n,k)$ to have  a non-expansive system $\mathcal{H}_{\infty}$ norm (i.e. $||G||_{\infty}\leq 1$). Before that, a specific arrangement of the reference vehicles in the platoon is introduced.
\begin{definition}
An arrangement of reference vehicles is called {\it{minimally dense}} (MD) if $\mathcal{P}(n,k)$ is partitioned into  line segments with length $2k+1$ starting from one end such that in the {\it {middle}} of each partition one reference vehicle is located (which will be connected to all of the followers in that partition). 
\end{definition}
Based on MD arrangement there exist $\lceil \frac{n}{2k+1} \rceil$ reference vehicles in $\mathcal{P}(n,k)$. The following theorem introduces conditions for $\mathcal{P}(n,k)$ to have non-expansive $\mathcal{H}_{\infty}$ norm. 

\begin{theorem}
Consider a $k$-nearest neighbor platoon $\mathcal{P}(n,k)$ with dynamics \eqref{eqn:errodrd}. If there exist at least $|\mathcal{R}|=\lceil \frac{n}{2k+1} \rceil$ reference vehicles, then there exists an arrangement of the reference vehicles satisfying $||G||_{\infty}\leq 1$. Moreover if the number of reference vehicles is less than $\lceil \frac{n}{2k+1} \rceil$, then there is no arrangement of reference vehicles satisfying $||G||_{\infty}\leq 1$.
\label{thm:suffff}
\end{theorem}
\begin{IEEEproof}
First, the sufficient condition is explored. Based on Corollary \ref{cor:corolllll}, a sufficient condition for $||G||_{\infty}\leq 1$ is to have  $\min_{i\in \mathcal{F}}\{\beta_i\}\geq 1$. By doing an MD arrangement of $\lceil \frac{n}{2k+1} \rceil$  reference vehicles in  $\mathcal{P}(n,k)$ we will have $\min_{i\in \mathcal{F}}\{\beta_i\}\geq 1$.

Next we have to show that with less than this number of reference vehicles, it is impossible to obtain $||G||_{\infty}\leq 1$. From a lower bound  in \eqref{eqn:hinfff}, a necessary condition for $||G||_{\infty}\leq 1$ is to have $\frac{|\mathcal{F}|}{|\partial \mathcal{R}|}\leq 1$. Based on the fact that $|\partial \mathcal{R}| \leq 2k|\mathcal{R}|$ we have $\frac{|\mathcal{F}|}{2k|\mathcal{R}|} \leq \frac{|\mathcal{F}|}{|\partial \mathcal{R}|}$. Thus $\frac{|\mathcal{F}|}{2k|\mathcal{R}|}\leq 1$ is a necessary condition for $||G||_{\infty}\leq 1$, which yields $|\mathcal{R}| \geq \frac{n}{2k+1}$.
\end{IEEEproof}
Based on Theorem \ref{thm:suffff}, the MD arrangement of reference vehicles in $\mathcal{P}(n,k)$ provides the minimum possible number of reference vehicles to yield a non-expansive $\mathcal{H}_{\infty}$ norm.

\subsection{Robustness to Time Delay}

Here we discuss the stability of  dynamics \eqref{eqn:errord} when the vehicles update their states  with a particular time delay $\tau \in [0,\tau_{max}]$. The dependency of the robustness of linear systems like \eqref{eqn:errord} to time delay is discussed in the following theorem, which is based on a general result in \cite{Buslowicz}.
 \begin{theorem}
The velocity tracking dynamics  \eqref{eqn:errord} is asymptotically  stable in the presence of constant time delay  $\tau  \in [0,\tau_{max}]$ if and only if
\begin{equation}
\tau_{\rm max}< \min_i\{\frac{\pi}{2\lambda_i(L_g)}\}=\frac{\pi}{2\lambda_{|\mathcal{F}|}(\mathcal{L}_g)}.
\label{eqn:deleq}
\end{equation}
\label{thm:delll}
\end{theorem}
Based on  Theorems \ref{thm:bounddmax} and \ref{thm:delll}, the following proposition introduces necessary and sufficient conditions for the stability of $\mathcal{P}(n,k)$ under time delay.
 \begin{proposition}
A vehicle platoon $\mathcal{P}(n,k)$ under velocity tracking dynamics \eqref{eqn:errord} in the presence of constant time delay $\tau  \in [0,\tau_{max}]$ is asymptotically stable if $\tau_{\rm max}\leq \frac{\pi}{8k}$ and it is unstable if  $\tau_{\rm max} > \frac{\pi}{2k}$.
\label{prop:chemidoonam}
\end{proposition}
\begin{IEEEproof}
According to Theorems \ref{thm:bounddmax} and  \ref{thm:delll}, necessary and sufficient conditions for asymptotic stability of \eqref{eqn:errord} in the presence of time delay are $\tau_{\rm max} < \frac{\pi}{2 d_{\rm max}^{\mathcal{F}}}$ and $\tau_{\rm max}  < \frac{\pi}{4 d_{\rm max}^{\mathcal{F}}}$, respectively, and according to the fact that $k \leq d^{\mathcal{F}}_{max} \leq 2k$ the results are obtained.
\end{IEEEproof}
\begin{remark}
As can be concluded from Corollary \ref{cor:corolllll} and Proposition \ref{prop:chemidoonam}, there is a trade-off between robustness to disturbances and time delay when the connectivity index $k$ increases. In particular, by increasing $k$ (for fixed number of reference vehicles) $\min_{i\in \mathcal{F}}\{\beta_i\}$ increases while the maximum delay $\tau_{\rm max}$ decreases.
\label{rem:tradeofrem1}
\end{remark}

\section{Robustness of Network Formation Dynamics}
\label{sec:doubleeee}

Similarly to Section \ref{sec:singleeee}, the performance and robustness of the network formation dynamics dynamics \eqref{eqn:doub} are analyzed in this section. 

\subsection{Stability Margin and Robustness to Uncertainty}

In a $k$-nearest neighbor platoon $\mathcal{P}(n,k)$, for the smallest magnitude of the real part of $\mathcal{B}$ (stability margin of \eqref{eqn:doub}) the following proposition is presented.

\begin{proposition}
For the stability margin of the network formation dynamics \eqref{eqn:doub} for the platoon $\mathcal{P}(n,k)$  we have
\begin{align}
|Re(\lambda_{1}(\mathcal{B}))| \geq \frac{\lambda_{1}(\mathcal{L}_g)}{2}.
\label{eqn:doubsmal}
\end{align}
\label{prop:boundsdelrela}
\end{proposition}
\begin{IEEEproof}
Based on Theorem \ref{thm:eigenthm} and \eqref{eqn:charact}, for the eigenvalue of $\mathcal{B}$ with the smallest magnitude of the real part two cases may occur: (i) If $1-\frac{4}{\lambda_{1}(\mathcal{L}_g)}\leq 0$ then $|Re(\lambda_{1}(\mathcal{B}))|=\frac{\lambda_{1}(\mathcal{L}_g)}{2}$, (ii) If $1-\frac{4}{\lambda_{1}(\mathcal{L}_g)}>0$ then $|Re(\lambda_{1}(\mathcal{B}))| > \frac{\lambda_{1}(\mathcal{L}_g)}{2}$.
\end{IEEEproof}

The following theorem gives an upper bound for the  $\mathcal{H}_{\infty}$ norm of the network formation dynamics  under the MD arrangement of reference vehicles.

\begin{theorem}
Consider $\mathcal{P}(n,k)$ with a MD arrangement of reference vehicles. The $\mathcal{H}_{\infty}$ norm from disturbances to the position error of  \eqref{eqn:doub} satisfies $||G||_{\infty} \leq \frac{2}{\sqrt{3}}$.
\label{thm:sqrt2}
\end{theorem}

\begin{IEEEproof}
Taking Laplace transform of \eqref{eqn:singaslle} for zero initial conditions gives
\begin{align}
G(s)=\left(s^2I+(s+1)\mathcal{L}_g \right)^{-1}
=M\left(s^2I+(s+1)\Lambda\right)^{-1}M^T =M \diag(G_i(s))M^T,
\end{align}
where  $M=[v_1, v_2, ..., v_{|\mathcal{F}|}]$ is a matrix formed by eigenvectors of $\mathcal{L}_g$ and $\diag(G_i(s))$ is a diagonal matrix with diagonal elements $G_i(s)=\frac{1}{s^2+\lambda_i(\mathcal{L}_g)s+\lambda_i(\mathcal{L}_g)}$  with the following maximum amplitudes 
\begin{equation}
\mathcal{C}_i=\max_{\omega}|G_i(j\omega)|=
  \begin{cases}
  \frac{2}{\lambda_i(\mathcal{L}_g)^{\frac{3}{2}}\sqrt{4-\lambda_i(\mathcal{L}_g)}},  & \quad \text{if } \lambda_i\leq 2,\\
   \frac{1}{\lambda_i(\mathcal{L}_g)} & \quad  \text{otherwise}.\\
  \end{cases}
  \label{eqn:casess}
\end{equation}
Hence, for system $\mathcal{H}_{\infty}$ norm we have
\begin{align}
||G||_{\infty}=\sup_{\omega}\max_i||G_i(j\omega)||=\max_i\mathcal{C}_i=\mathcal{C}_1.
\end{align}
Now due to the fact that in MD arrangement we have $1 \leq \lambda_1(\mathcal{L}_g)\leq 2$, and considering the fact that in this interval $\mathcal{C}_1$ takes its maximum at $\lambda_1(\mathcal{L}_g)=1$ we have
\begin{align}
 ||G||_{\infty}=\frac{2}{\lambda_1(\mathcal{L}_g)^{\frac{3}{2}}\sqrt{4-\lambda_1(\mathcal{L}_g)}} \leq \frac{2}{\sqrt{3}}.
\end{align}
\end{IEEEproof}

Theorems \ref{thm:suffff} and \ref{thm:sqrt2} show how the existence of  multiple reference vehicles in a platoon can increase the robustness of the network against disturbances. In Table \ref{tab:nscvh}, the system $\mathcal{H}_{\infty}$ norm of the velocity tracking and network formation dynamics on $\mathcal{P}(n,k)$ for both single and multiple reference vehicles with MD arrangement, i.e. $|\mathcal{R}|=\lceil \frac{n}{2k+1} \rceil$, is summarized.\footnote{In \cite{hao2013stability} it is shown that the  system $\mathcal{H}_{\infty}$ norm for network formation dynamics for $\mathcal{P}(n,1)$ (line graph) is $\Theta({n^3})$, which holds for any $k<\infty$ as well. Moreover, it can be easily shown that  the $\mathcal{H}_{\infty}$ norm of the velocity tracking dynamics is $\Theta({n^2})$, due to the fact that for line graphs we have $\lambda_1(\mathcal{L}_g)=\Theta(\frac{1}{n^2})$ \cite{Yueh}.}

\begin{table}[h] 
\centering
\caption {System $\mathcal{H}_{\infty}$ norm of \eqref{eqn:errord} and \eqref{eqn:doub}  for single and multiple reference vehicles with MD arrangement.}
\begin{tabular}{c c c c c} 
 \hline
 $|\mathcal{R}|$ & Velocity Tracking \eqref{eqn:errord} &  Network Formation \eqref{eqn:doub} \\
 \hline
 $1$  & $\Theta({n^2})$ & $\Theta({n^3})$ \\
 $\lceil \frac{n}{2k+1}\rceil$ & $\leq 1$ & $\leq \frac{2}{\sqrt{3}}$ \\
 \hline\\
\end{tabular}
\label{tab:nscvh}
\end{table}

\subsection{Robustness to Time Delay}

The following proposition gives a sufficient condition for which the network formation dynamics \eqref{eqn:doub} remains asymptotically stable in the presence of time delay. 

\begin{proposition}
The network formation  dynamics  \eqref{eqn:doub} is  asymptotically stable in the presence of constant time delay  $\tau  \in [0,\tau_{\rm max}]$  if $\tau_{\rm max}<\frac{1}{4k}$. 
 \label{prop:delcor}
\end{proposition}
 \begin{IEEEproof}
 Based on \cite{Niculescu}, a sufficient condition for \eqref{eqn:doub} to remain stable in the presence of time delay is to have
\begin{equation}
\tau_{\rm max}<\frac{1}{\rho(\mathcal{B})},
\label{eqn:sufrhodoub}
\end{equation}
where $\rho(\mathcal{B})$ is the spectral radius of $\mathcal{B}$. Applying Theorem \ref{thm:eigenthm} and \eqref{eqn:charact}, the spectral radius  of $\mathcal{B}$ yields: 
 \begin{equation}
 \rho(\mathcal{B})=\frac{\lambda_{|\mathcal{F}|}(\mathcal{L}_g)}{2}\left(1+\left(1-\frac{4}{\lambda_{|\mathcal{F}|}(\mathcal{L}_g)}\right)^{\frac{1}{2}} \right),
 \end{equation}
since in the MD arrangement in which $\lambda_1(\mathcal{L}_g)\geq 1$ we have $\max_i |1-\frac{4}{\lambda_{i}(\mathcal{L}_g)}| = |1-\frac{4}{\lambda_{|\mathcal{F}|}(\mathcal{L}_g)}|$. Therefore, based on the upper bound on $\lambda_{|\mathcal{F}|}(\mathcal{L}_g)$ in \eqref{eqn:up}, sufficient condition \eqref{eqn:sufrhodoub} can be rewritten as
 \begin{equation}
 \tau_{\rm max}<\frac{1}{d_{\rm max}^{\mathcal{F}} +d_{\rm max}^{\mathcal{F}}\left(1-\frac{2}{d_{\rm max}^{\mathcal{F}} }\right)^{\frac{1}{2}}}
 \end{equation}
and since $d_{\rm max}^{\mathcal{F}} \geq 2$ we have $d_{\rm max}^{\mathcal{F}} +d_{\rm max}^{\mathcal{F}}\left(1-\frac{2}{d_{\rm max}^{\mathcal{F}} }\right)^{\frac{1}{2}}\leq 2d_{\rm max}^{\mathcal{F}}$. This yields the sufficient condition $ \tau_{\rm max}<\frac{1}{2d_{\rm max}^{\mathcal{F}}}$,
and based on the fact that $d_{\rm max}^{\mathcal{F}}\leq 2k$ the result will be obtained. 
 \end{IEEEproof}
 \begin{remark}
 Similar to what was mentioned in Remark \ref{rem:tradeofrem1} for velocity tracking dynamics, there is a trade-off between robustness to disturbances and time delay for the network formation dynamics. More specifically, by increasing network connectivity $k$ the value of $\lambda_1(\mathcal{L}_g)$ increases and based on \eqref{eqn:casess} the system $\mathcal{H}_{\infty}$ norm decreases, while  the spectral radius of  $\mathcal{B}$  increases which results in decreasing the robustness to time delay.  
 \end{remark}

\section{Simulations} \label{sec:Results}
In this section, some simulation results are presented to confirm the theoretical contributions of the paper. The results are based on  $\mathcal{P}(36,4)$. Based on MD arrangement, there are four reference vehicles in $\mathcal{P}(36,4)$.

\begin{figure}[h!]
\centering
\includegraphics[width=0.6\linewidth]{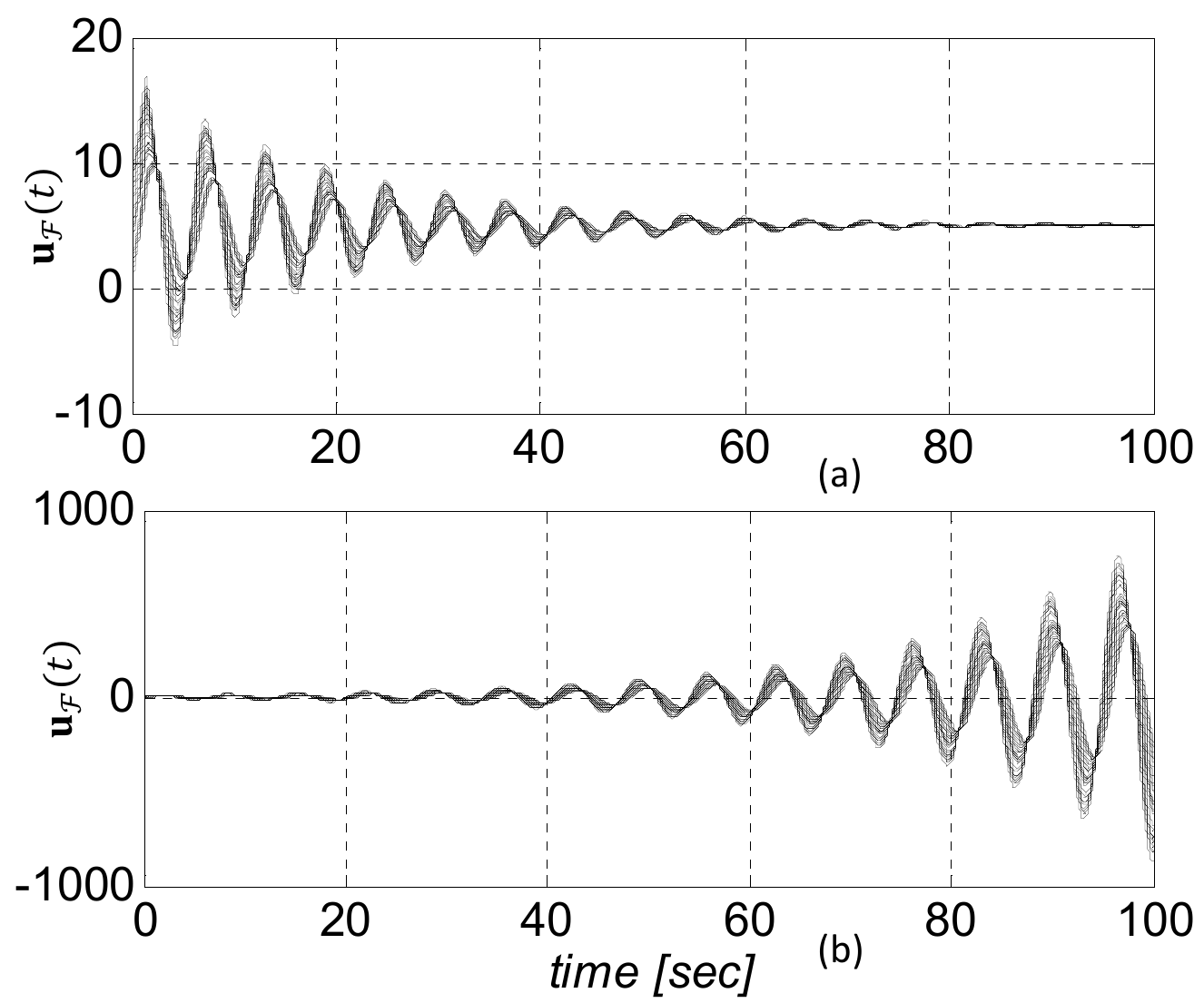}
\caption{(a) velocity tracking dynamics is stable for $\tau=0.09<\frac{\pi}{2k}$ and (b) unstable for $\tau=0.4>\frac{\pi}{8k}$.}
\label{fig:SI_n36_k4_Delay_Stable_Unstable}
\end{figure}

\begin{figure}[h!]
\centering
\includegraphics[width=0.6\linewidth]{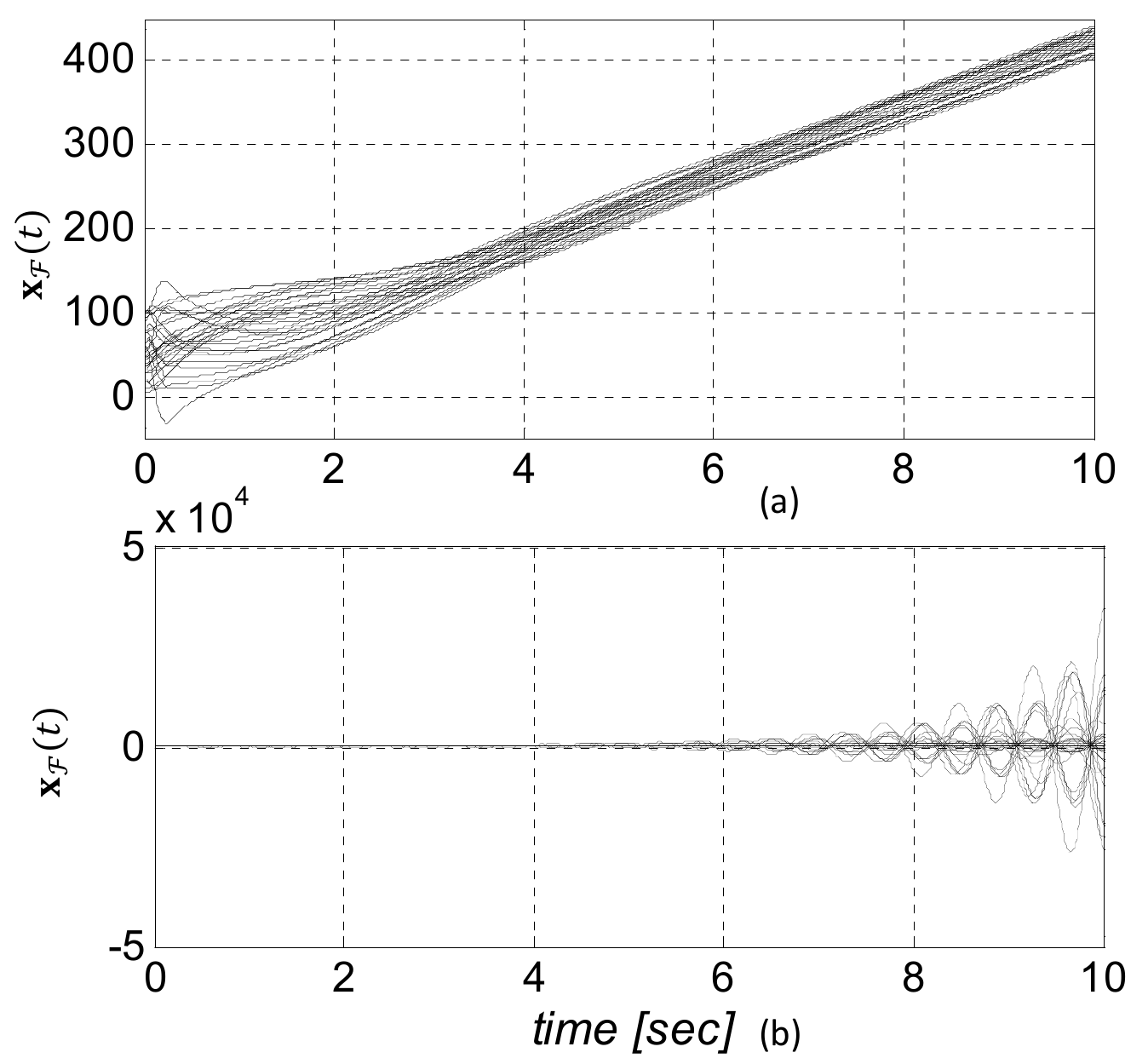}
\caption{(a) network formation dynamics is stable for $\tau=0.05<\frac{1}{4k}$ and (b) unstable for $\tau=0.1>\frac{1}{4k}$.}
\label{fig:DI_n36_k4_Delay_Stable_Unstable}
\end{figure}

\begin{figure}[h!]
\centering
\includegraphics[width=0.6\linewidth]{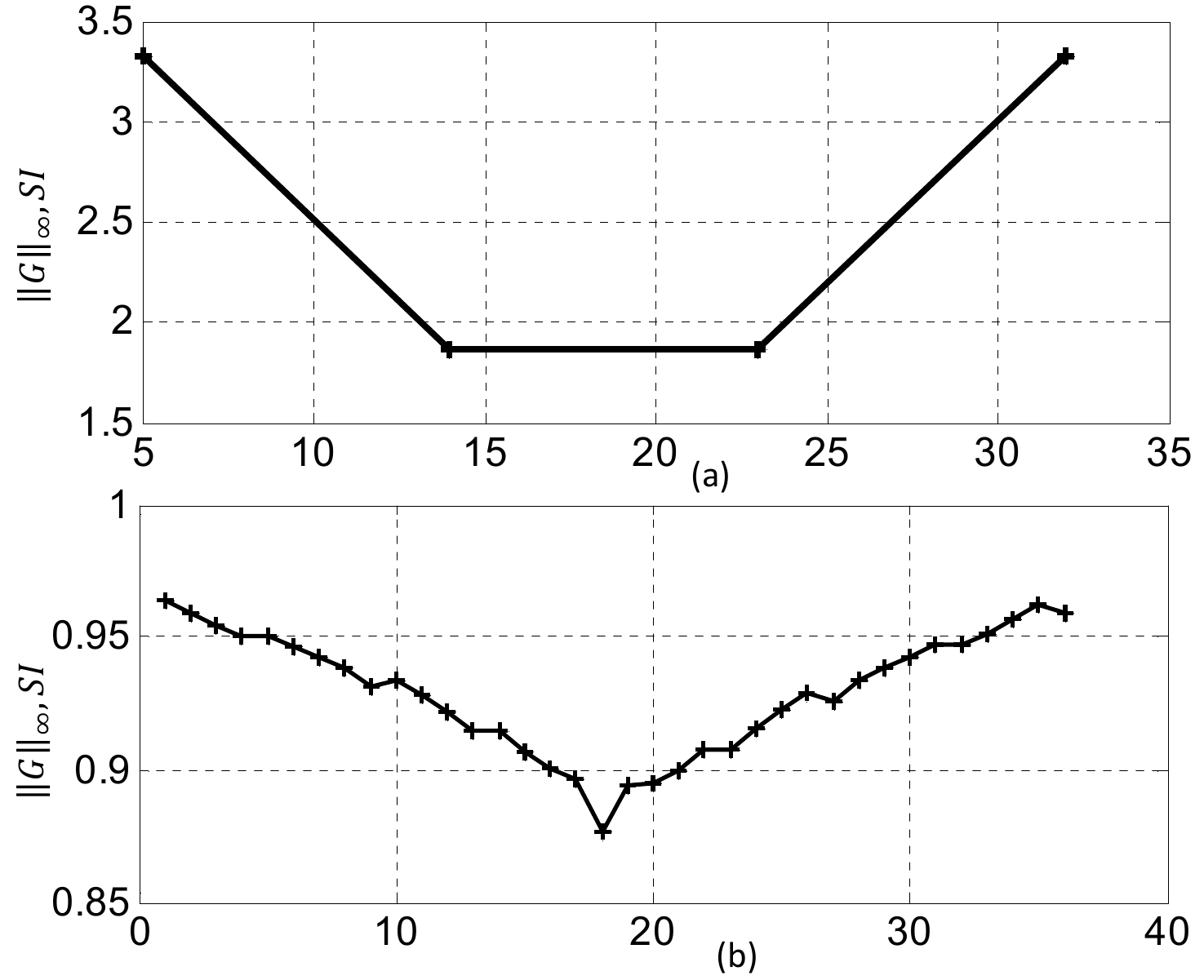}
\caption{(a) $\mathcal{H}_{\infty}$ norm of \eqref{eqn:errord} with removing a single reference vehicle from the MD arrangement in $\mathcal{P}(36,4)$ and (b) adding a single reference vehicle. The horizontal axis is the location of the added (removed) vehicle in the platoon.}
\label{fig:SI_Hinfn36_k4_Nec_Suf}
\end{figure}

\begin{figure}[h!]
\centering
\includegraphics[width=0.6\linewidth]{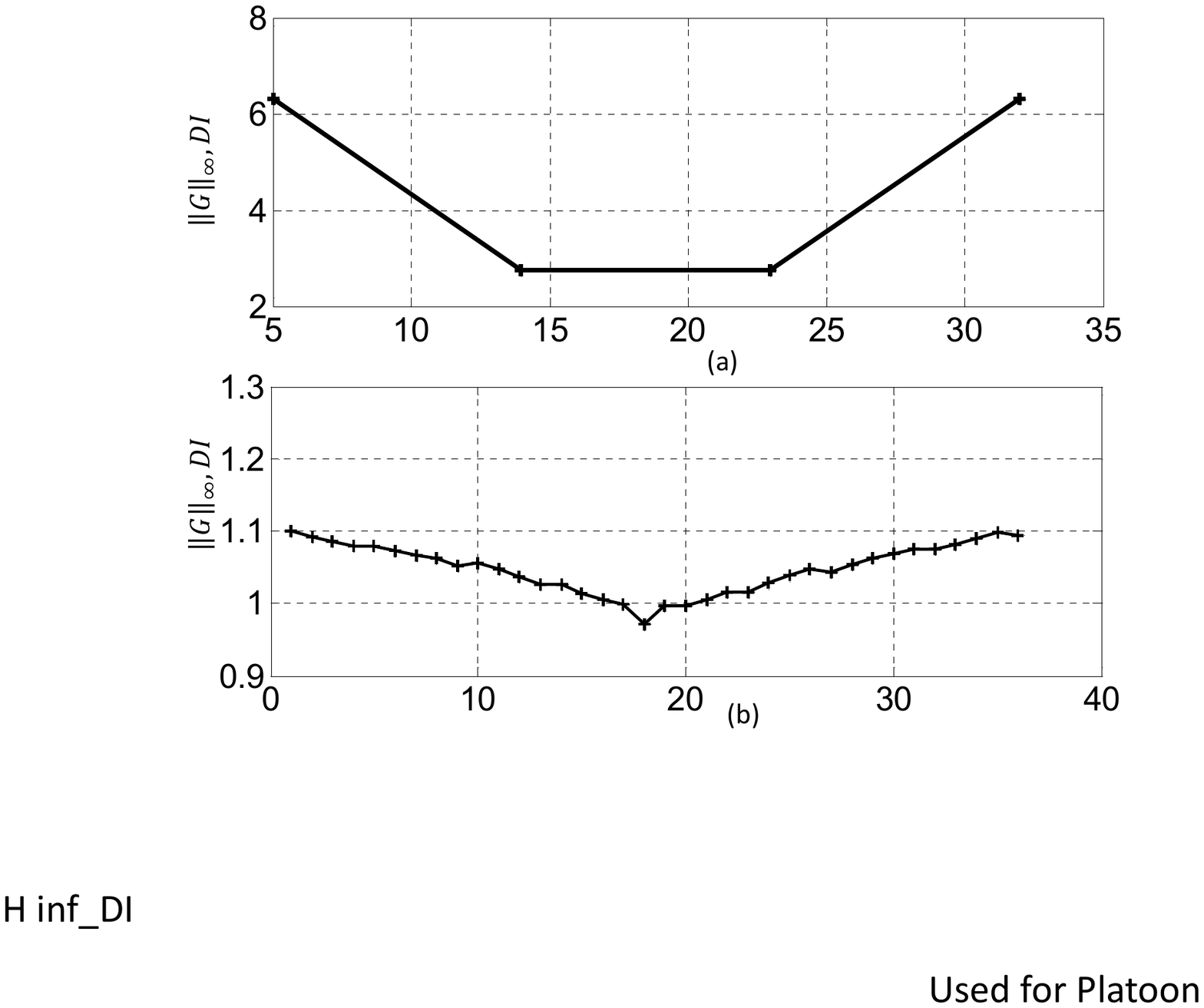}
\caption{(a) $\mathcal{H}_{\infty}$ norm of \eqref{eqn:doub} with removing a single reference vehicle from the MD arrangement in $\mathcal{P}(36,4)$ and (b) adding a single reference vehicle and its location
.}
\label{fig:DI_Hinfn36_k4_Nec_Suf}
\end{figure}

Fig. \ref{fig:SI_n36_k4_Delay_Stable_Unstable} shows how necessary and sufficient conditions for the value of time delay mentioned in Proposition \ref{prop:chemidoonam} apply for asymptotic stability of  the velocity tracking dynamics  \eqref{eqn:errord} in the presence of time delay. The sufficient condition for the stability of the network formation dynamics \eqref{eqn:doub} in the presence of time delay (presented in Proposition \ref{prop:delcor}) is confirmed in Fig.~\ref{fig:DI_n36_k4_Delay_Stable_Unstable}. 
Fig.~\ref{fig:SI_Hinfn36_k4_Nec_Suf} shows how MD arrangement of the reference vehicles introduces tight necessary and sufficient conditions for system $\mathcal{H}_{\infty}$ norm of \eqref{eqn:errord} to be non-expansive (Theorem \ref{thm:suffff}). In particular, if one of the four reference vehicles in the MD arrangement of $\mathcal{P}(n,k)$ is removed, the resulting $\mathcal{H}_{\infty}$ norm is no longer less than one. On the other hand, as can be seen from Fig.~\ref{fig:SI_Hinfn36_k4_Nec_Suf}  if an extra reference vehicle is added to $\mathcal{P}(n,k)$ (other than the existing reference vehicles from the MD arrangement), the resulting $\mathcal{H}_{\infty}$ norm will be strictly less than one. The results for the same scenario are shown for the network formation dynamics \eqref{eqn:doub} as shown in Fig.~\ref{fig:DI_Hinfn36_k4_Nec_Suf} where removing a reference vehicle makes the $\mathcal{H}_{\infty}$ norm of \eqref{eqn:doub} larger than $\frac{2}{\sqrt{3}}\approx 1.15$. This confirms the tight sufficient condition mentioned in Theorem \ref{thm:sqrt2}.

\section{Summary and Conclusions}
 \label{sec:conclusion}
 In this paper a set of graph theoretic conditions for the robustness of $k$-nearest neighbor vehicle platoons $\mathcal{P}(n,k)$ to disturbances and time delay have been derived and analyzed. In particular, a necessary and sufficient condition for $\mathcal{P}(n,k)$ to have non-expansive $\mathcal{H}_{\infty}$ norm for velocity tracking dynamics has been provided (Theorem \ref{thm:suffff}) by introducing a specific arrangement of reference vehicles. Moreover, the effect  of such arrangement of reference vehicles on $\mathcal{H}_{\infty}$ norm of network formation dynamics has been investigated (Theorem \ref{thm:sqrt2}). Furthermore, the effect of the communication delay on the stability of velocity tracking dynamics and network formation dynamics has been addressed (Propositions \ref{prop:chemidoonam} and \ref{prop:delcor}). Our results show that there is a trade-off between robustness to time delay and robustness to disturbances. An avenue for future work in this direction is to generalize the results established in this paper to directed networks with non-homogeneous control gains.


\bibliographystyle{IEEEtran}
\bibliography{main.bib}

\end{document}